\documentclass[1pt,notitlepage,twoside,a4paper]{amsart}

\usepackage{amsmath,amssymb,enumerate}

\usepackage{epsfig,fancyhdr,color}

\usepackage{amssymb}
\usepackage{amsmath,amsthm}    
\usepackage{latexsym}
\usepackage{amscd} 
\usepackage{psfrag}
\usepackage{graphicx} 
\usepackage[latin1]{inputenc}  
\usepackage[all]{xy} 
\usepackage[mathcal]{eucal}

\theoremstyle{definition}

\theoremstyle{remark}

\def\interieur#1{\mathord{\mathop{\kern 0pt #1}\limits^\circ}}

\definecolor{NoteColor}{rgb}{1,0,0}

\title[A Commentary on Teichm\"uller's paper]{A Commentary on Teichm\"uller's paper \emph{Vollst\"andige L\"osung einer Extremalaufgabe der quasikonformen Abbildung (Complete solution of an extremal problem of the quasiconformal mapping)}}

\author[Vincent Alberge and Athanase Papadopoulos]{Vincent Alberge and Athanase Papadopoulos}

\address{Institut de Recherche Math\'ematique Avanc\'ee\\ CNRS et Universit\'e de Strasbourg\\\small 7 rue Ren\'e
  Descartes - 67084 Strasbourg Cedex, France;
  \\
  The Graduate Center, City University of New York,\\
365 Fifth Avenue,
New York, NY 10016 USA. \\
\tt{alberge@math.unistra.fr}, \tt{papadop@math.unistra.fr}}

 \date{\today}

\begin{document}

 \begin{abstract} 

We comment on Teichm\"uller's paper \emph{Vollst\"andige L\"osung einer Extremalaufgabe der quasikonformen Abbildung (Complete solution of an extremal problem of the quasiconformal mapping)} \cite{T24}, published in 1941. In this paper, Teichm\"uller gives a proof of the existence of extremal quasiconformal mappings in the case of the pentagon (disc with five distinguished points on the boundary).

 \medskip

 \medskip

 \medskip

The final version of this paper will appear as a chapter in Volume VI of {\it the Handbook of Teichm\"uller theory}. This volume is dedicated to the memory of Alexander Grothendieck.

\end{abstract}

\bigskip

\maketitle

\tableofcontents

\begin{figure}[ht!]
\centering
 \psfrag{0}{\small $0$}
 \psfrag{a}{\small $a$}
 \psfrag{-b}{\small $-b$}
 \psfrag{A}{\small $A$}
\includegraphics[width=.80\linewidth]{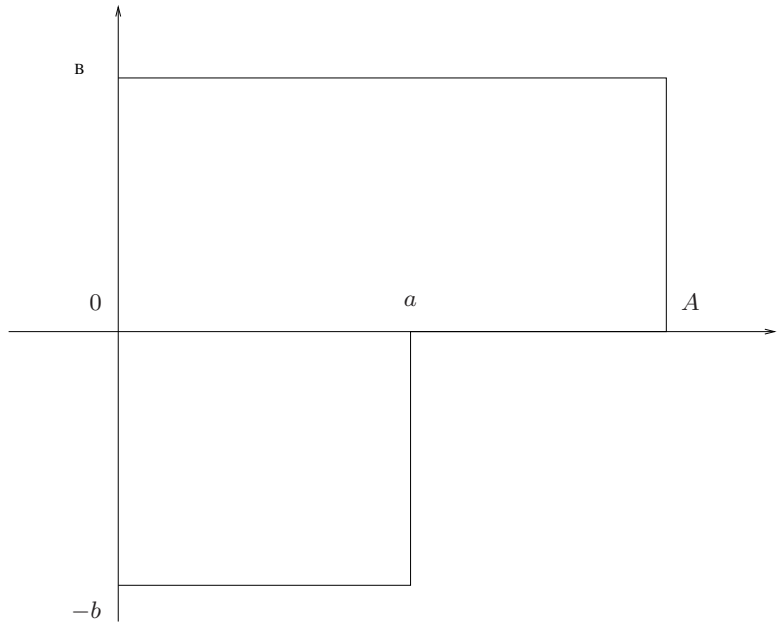}
\caption{\small{This hexagonal-shaped figure is a pentagon (the five distinguished points are at the salient angles).}}
\label{hexagon2}
\end{figure}

We comment on the paper \emph{Vollst\"andige L\"osung einer Extremalaufgabe der quasikonformen Abbildung (Complete solution of an extremal problem of the quasiconformal mapping)} \cite{T24} by Teichm\"uller, published in 1941. In this paper Teichm\"uller proves his famous \emph{existence theorem} of extremal quasiconformal mappings (we shall recall the statement below), for the case of the pentagon (a disc with five distinguished points on the boundary). In the paper \cite{T20} (see also the commentary \cite{T20C}), published in 1939  and which is probably his most quoted paper, Teichm\"uller had announced this theorem for arbitrary surfaces of finite topological type (orientable or not), with only a sketch of a proof. A complete proof of this theorem is given in that paper only in the case of the torus and in a few other cases which can be reduced to that case: the sphere with four distinguished points, the annulus (a case which was already treated by Gr\"otzsch), the disc with two interior distinguished points, and the disc with one interior distinguished point and two boundary distinguished points.
 In the later paper \cite{T29} (see also the commentary \cite{T29C}), published in 1943, Teichm\"uller gave a complete proof of this theorem in the case of an orientable closed surfaces of finite type. The case of the pentagon is not proved rigorously in  the paper \cite{T20} (and this case is also not considered in the paper \cite{T29} which was published later).\footnote{In \cite{T29}, Teichm\"uller promises to give later on a proof in the most general case of surfaces of finite type (orientable or not, with or without boundary, with or without distinguished points in the interior and/or on the boundary). His project was not realized since he died soon later.} However, Teichm\"uller proves in \S\,129 of  \cite{T20}, using the length-area method, that the affine maps (which are in fact Teichm\"uller maps) between two pentagons of a particular type (see Figure \ref{hexagon2}) are extremal. Furthermore, in \S\,130 of the same paper, he develops the geometry of the space we call today the Teichm\"uller space of the pentagon,\index{pentagon!Teichm\"uller space}\index{Teichm\"uller space!pentagon} equipped with its Teichm\"uller metric. He studies quasiconformal maps between pentagons and the geodesics in this Teichm\"uller space. The existence theorem for extremal quasiconformal mappings between two pentagons is implicitly admitted (and not proved) in that paper. This is why the present paper is a valuable addition to the paper  \cite{T20}. It is also important to add that in \cite{T20}, Teichm\"uller was still not using the so-called method of continuity,\index{method of continuity} and this is why in \S\,161 to \S\,163 of that paper he tried to use the so-called length-area method\index{length-area method}\footnote{In \cite{T20}, Teichm\"uller calls this method the Gr\"otzsch-Ahlfors method.\index{Gr\"otzsch-Ahlfors method}}  in order to determine the extremal map. His idea was probably to generalize the method used for the quadrilateral which originates in the work of Gr\"otzsch.

Let us note by the way that the case of the pentagon is a nontrivial one. Gr\"otzsch, in 1932, had treated the case of a quadrilateral (a disc with four distinguished points on the boundary)  \cite{Gr1932}. In his 1964 survey paper on quasiconformal maps and their applications \cite{Ahlfors1964}, Ahlfors, reporting on Teichm\"uller's work,  writes that the result on pentagons is ``already a sophisticated result.'' Teichm\"uller writes (\S 1 of the present paper) that the case of the pentagon is ``the simplest case of the higher cases,'' and that this simple case already ``shows how far one has to go beyond and extend the methods of Ahlfors and Gr\"otzsch.'' Indeed, some nontrivial work is needed  for the existence proof in this case, and the proofs  already highlight the difficulties that appear in the general case. In fact, Teichm\"uller wrote explicitly that he will only deal with the ``solution of the problem of extremal quasiconformal mapping in the case of the pentagon,'' and that to solve this problem he will use a proof by continuity; a method which ``may serve [...] for a proof of the general case.'' We know that Teichm\"uller used  this method in \cite{T29} to give a rigorous proof for the existence of extremal quasiconformal mappings in the case of  closed surfaces of genus $\geq 2$.  We also note that the Teichm\"uller space of the pentagon,  which is studied in some detail in the paper \cite{T20}, coincides with the moduli space of this surface since the mapping class group in this case is trivial (the distinguished points on the boundary are pointwise fixed by the mappings).

Let us make a further remark on the geometry of pentagons.

One of the beautiful results of Teichm\"uller's work is that 
the extremal quasiconformal mappings (the so-called Teichm\"uller mappings) between arbitrary surfaces of finite type are locally affine; more precisely, in the local $\zeta$-coordinate of the complex plane, they have the form 
\[
\zeta \to K\cdot\mathrm{Re}\left( \zeta \right) + \mathrm{Im}\left( \zeta \right).\]
In the case of the quadrilateral, this is also the form of a global map, by a result of Gr\"otzsch, obtained after conformally mapping an arbitray quadrilateral to a Euclidean rectangle with sides parallel to the real and imaginary axes of the complex plane. Now if one hopes for such a result for more general surfaces, one needs to find good conformal representatives of more general polygons (discs with distinguished points on the boundary). Teichm\"uller succeeds in doing this, but representing a pentagon by a Euclidean \emph{hexagon}, that is, a figure with 6 vertices instead of five (see Figure \ref{hexagon2} below). In this case, one of the vertices (the one with the re-entering angle) is not considered as a distinguished point, and the hexagon becomes the conformal image of a pentagon.

Let us now present in some detail the results in this paper.

In \S 2 of the paper, Teichm\"uller states precisely the existence problem for extremal quasiconformal mappings, in the case of pentagons. A pentagon is a disc with five ordered distinguished points on its boundary.  It is represented conformally (using the Riemann mapping theorem) as the upper half-plane $\mathrm{Im}\left( z\right)>0$ together with its boundary, such that (using Teichm\"uller's notation) the distinguished points are sent to the ordered quintuple $0,p_2,1,p_4,\infty$. Thus, the pair of real numbers $p_2,p_4$, with the conditions $0<p_2<1$, $1<p_4 <\infty$, are the parameters for the space of conformal classes of pentagons. Teichm\"uller called such pentagons ``normalized pentagons.'' The quasiconformal mappings between pentagons are defined, as continuously differentiable mappings in both directions, with singularities at finitely many points and at finitely many analytical arcs. The \emph{dilatation quotient}\index{dilatation quotient} is defined as the ratio of the large axis to the small axis of the infinitesimal ellipse which is the image of an infinitesimal circle at a point where the mapping is differentiable, and the map is said to be quasiconformal if the supremum of the dilatation quotient over the whole surface (whenever the quotient is defined) is finite.
Teichm\"uller assumes the reader familiar with this notion, and he refers to his earlier paper \cite{Teine}. We call this supremum the \emph{quasiconformal dilatation}\index{quasiconformal dilatation} of the map.

Given two pentagons $P,Q$ there always exists a quasiconformal mapping between them, and Teichm\"uller gives the following example of such a mapping. The upper half-plane is sent conformally to the unit disc. The five distinguished points divide the boundary of this disc into five sectors. Using the polar coordinates $\rho e^{i\theta}$, the corresponding sectors are sent to each other by maps of the form $\rho'=\rho$ and $\theta'=a\theta+b$, with $a>0$. The problem is now the following: 
\begin{quote} Given two pentagons $P$ and $Q$, is there an \emph{extremal} quasiconformal mapping, that is, a mapping  which has the smallest possible quasiconformal dilatation?\end{quote}
The aim of this paper is to answer this question affirmatively. Furthermore, Teichm\"uller shows that the dilatation quotient of this extremal map is everywhere constant, and he gives an explicit expression for the form that this map has.

The proof will use the so-called ``continuity argument,'' in the following form. One starts with an arbitrary pentagon $P$, and builds a two-parameter family of pentagons $Q$, equipped with special mappings $P\to Q$. These special mappings are shown to be extremal quasiconformal, and the ``continuity argument'' will show that an arbitrary pentagon is obtained in this manner. The continuity argument is Brouwer's theorem of invariance of domain.\index{theorem!invariance of domain}

In \S 3, titled \emph{Maps to axis-parallel hexagons}, Teichm\"uller uses a  mapping given by the integral
\begin{equation}\label{schwarz}\zeta=\int \sqrt{\frac{\cos \varphi+z\sin \varphi}{z(z-p_2)(z-1)(z-p_4)}}dz
.\end{equation}
Such a mapping is called a \emph{Schwarz-\-Christoffel mapping}\index{Schwarz-Christoffel mapping}.
It sends the family of polygons parametrized as above by the points $p_2,p_4$ on the boundary of the upper half-plane to a family of polygons in the $\zeta$-plane which are simply connected regions having the form of Euclidean polygons. Depending on the values of the parameter $\varphi$, the image is either a Euclidean hexagon (when $\cot\left( \varphi\right)\neq 0, -p_2 , -1,  -p_4 , \infty$) where distinguished points correspond to vertices with angle $\frac{\pi}{2}$, this Euclidean hexagon having also a re-entrant vertex of angle $\frac{3\pi}{2}$ (this vertex is not considered as a distinguished point), or a rectangle with five marked points  where the four vertices are distinguished points, and one of the five distinguished points lies on one of the sides. See Figure \ref{hexagon2}. Using the notation in that figure, by making  $a=A$, or $b=0$, the hexagon degenerates into a rectangle. Teichm\"uller calls this kind of degenerate hexagon a  \emph{rectangular hexagon}. We recall that the study of pentagons using a representation by Euclidean hexagons was already done in \S\,129 of \cite{T20}. For a fixed pentagon determined by $p_2,p_4$, taking the square of the above integrand,
\[d\zeta^2= \frac{\cos \varphi+z\sin \varphi}{z(z-p_2)(z-1)(z-p_4)}dz^2
,\]
is, up to a positive factor, the general form of the real and not identically zero quadratic differentials on the Riemann sphere represented by the parameter $z$ and the point at infinity. This kind of quadratic differential is real along the real axis and is meromorphic with at most simple poles at the five points $0, p_2 , 1, p_4 , \infty$. Teichm\"uller calls such a quadratic differential ``regular."\footnote{The explicit definition of such a notion is given in \S\,100 of \cite{T20}.}

In \S\,4, a space of pairs $(P,\varphi)$ is studied. It has three parameters: two parameters for $(p_2,p_4)$  varying in $(0,1)\times (1,\infty)$, and one for the angle parameter $\varphi$. This spaces is equipped with the canonical topology induced by that of $\mathbb{R}^3$. A map $(P,\varphi)\to S$ is defined from this parameter space into the space  of equivalence classes of hexagons (and rectangular hexagons) with the parameters $a,b, A,B$ that we mentioned, up to the transformations
$$
\zeta \to a\zeta + b,
$$
where $b\in\mathbb{C}$ and $a\in\mathbb{R}^*$. The topology on this space of pentagons (which have the shape of hexagons) is the one induced by the 
four numbers $a,b, A,B$ up to a scalar factor, or, as Teichm\"uller puts it, by the homogeneous coordinates $a:b:A:B$.  He makes a detailed study of the cases where the hexagonal shapes degenerate to rectangles.  
Teichm\"uller shows that the map $(P,\varphi)\to S$  is continuous. For that, he starts by showing that the image of the map is a 3-dimensional topological manifold.  Using the dominated convergence theorem for integrals, he deduces that the map given by the integral of type 
 (\ref{schwarz}) depends continuously on $\left( p_2 , p_4 , \varphi \right)$. The next goal is to show that this map is a homeomorphism. This is done in the next section.

 In \S\,5, Teichm\"uller starts by showing that the map given by the integral is injective. Thus, this map is continuous and one-to-one onto its image. He then applies Brouwer's theorem of invariance of domain to conclude \index{theorem!Invariance of domain} that the map, in fact, is a homeomorphism onto its image. Teichm\"uller notes however that the use of this theorem ``can still be avoided.'' He concludes this section by showing that the map is in fact surjective. Therefore the map $(P,\varphi)\to S$ is a homeomorphism.
 
 In \S\,6, Teichm\"uller defines what is called today the \emph{Teichm\"uller mapping}. Let us recall the construction.  Let $P$ be a pentagon determined by a pair $\left( p_2 , p_4 \right)\in (0,1)\times (1,\infty)$. Using the formula given by (\ref{schwarz}), we obtain, from $\varphi$ (a real number modulo $2\pi$), a new coordinate   $\zeta$ in which the pentagon has one of the desired forms, that is, either a hexagon in the plane whose distinguished points are the five salient vertices,  or a rectangular hexagon. In either case, we obtain an equivalence class of such hexagons called $S$. Teichm\"uller defines, for $K\geq 1$, a map (which is called now a \emph{Teichm\"uller map}\index{Teichm\"uller map!pentagon}\index{pentagon!Teichm\"uller map}), between two Euclidean hexagons which in natural local coordinates has the form
\begin{equation}\label{teichmap}
\zeta \to K\cdot\\mathrm{Re}\left( \zeta \right) + \mathrm{Im}\left( \zeta \right).
\end{equation}
This map is quasiconformal and its dilatation quotient is equal to $K$. Via the homeomorphism that we already considered, we get a new pentagon $Q$, denoted by $P\left( K, \varphi\right)$. Thus, Teichm\"uller defines a map
\begin{equation}\label{teichtheo}
\left( K, \varphi \right) \to P\left( K, \varphi\right)
\end{equation}
with values in the Teichm\"uller space of the pentagon. We have 
\begin{equation*}
P\left( 1, \varphi \right)= P. 
\end{equation*}

In the next section, Teichm\"uller shows that the map given by   (\ref{teichtheo}) is continuous. This is only a consequence of the fact that the set of pairs $\left( P, \varphi \right)$ is homeomorphic to the set of equivalence classes of Euclidean hexagon-shaped figures. 

In \S\,8, Teichm\"uller shows that for a given $\left(K, \varphi\right)$, the quasiconformal mappings between $P$ and $P\left( K, \varphi\right)$ that are induced from (\ref{teichmap}) are extremal. This result (in its general form, for an arbitrary surface of finite type) is often called the \emph{Teichm\"uller uniqueness theorem}.\index{Teichm\"uller uniqueness theorem}\index{Theorem!Teichm\"uller uniqueness!pentagon}\index{pentagon!Teichm\"uller uniqueness theorem}
 For the proof, Teichm\"uller uses the so-called length-area method.\index{length-area method} As we already noted, this was already proved in \S\,129 of \cite{T20} for pentagons and in \S\,132 to \S\,140 of that paper in full generality. The novelty in the present paper is in the next section.

Indeed, in the last section, Teichm\"uller shows that the map (\ref{teichtheo}) is a homeomorphism. To do this, he starts by showing that it is injective. Given that the map is continuous (\S\,7), he concludes, again using Brouwer's  theorem of invariance of domain\index{theorem!Invariance of domain} that this map is a homeomorphism onto its image. He then shows that the image corresponds to $\left\lbrace \left( p_2 , p_4 \right) \, \mid 0<p_2 <1 , 1< p_4 <\infty \right\rbrace$, i.e. the Teichm\"uller space of the pentagon. This concludes the \emph{Teichm\"uller theorem}\index{Teichm\"uller theorem!pentagon}\index{theorem!Teichm\"uller!pentagon}\index{pentagon!Teichm\"uller theorem}
 for the pentagon. 

As a conclusion to this commentary, let us note that the proof of the Teichm\"uller theorem for general compact hyperbolic surfaces given in \cite{T29} is modelled on the same idea, namely the application of Brouwer's theorem of invariance of domain to a space which is homeomorphic to Teichm\"uller space.

In the paper \cite{Ahlfors1946}, Ahlfors and Beurling consider conformal parameters for the pentagon and for some other special planar surfaces. This was the beginning of the notion of extremal length,\index{extremal length} which they try to apply to the study of the moduli of discs with a few number of distinguished points on the boundary, after representing them as polygons in the Euclidean plane.  Ahlfors and Beurling do not quote the paper \cite{T24} (it is possible that they were not aware of it), but they quote \cite{T20}. They write the following: ``Those who are familiar with the beautiful works of Mr. Teichm\"uller will notice the link between our results and his results. Nevertheless, one has to note that our starting point and the problem in which we are interested are enough far from the notions with which Mr. Teichm\"uller was concerned."\footnote{``Ceux qui connaissent bien les beaux travaux de M. Teichm\"uller noterons [sic] le lien \'etroit entre nos r\'esultats et les siens. Il faut tout-de-m\^eme remarquer que notre point de d\'epart ainsi que le probl\`eme qui nous int\'eresse sont assez \'eloign\'es des notions dont M. Teichm\"uller s'est pr\'eoccup\'e."}

\medskip

 \noindent {\bf Acknowledgements} The authors were partially supported by the French ANR program FINSLER and by the U.S. National Foundation grants DMS 1107452, 1107263, 1107367 ``RNMS  GEometric structures And Representation varieties.'' They wish to thank the Graduate Center of the City University of New York  where part of this work was done.

\end{document}